\documentclass[a4paper,12pt]{article}
\usepackage{latexsym}
\usepackage{amsmath}
\usepackage{amssymb}
\usepackage{enumerate}

\usepackage{theorem}
\newtheorem{theorem}{Theorem}[section]

\newtheorem{question}{Question}

\theorembodyfont{\rmfamily}

\newtheorem{remark}[theorem]{Remark}

\makeatletter

\@addtoreset{equation}{section}
\makeatother

\newcommand{\qedd}{\hfill \Box}
\newcommand{\ve}{\varepsilon}

\newcommand{\lra}{\longrightarrow}
\newcommand{\del}{\partial}

\newcommand{\R}{\ensuremath{\mathbb{R}}}

\newcommand{\bS}{\ensuremath{\mathbf{S}}}
\newcommand{\CD}{\mathsf{CD}}

\def\vol{\mathop{\mathrm{vol}}\nolimits}

\def\Ric{\mathop{\mathrm{Ric}}\nolimits}

\title{Vanishing S-curvature of Randers spaces}
\author{Shin-ichi OHTA\thanks{
Supported in part by the Grant-in-Aid for Young Scientists (B) 20740036.}\\
{\normalsize Department of Mathematics, Faculty of Science, Kyoto University,}\\
{\normalsize Kyoto 606-8502, JAPAN (e-mail: {\sf sohta@math.kyoto-u.ac.jp})}}
\begin{document}

\maketitle

\begin{abstract}
We give a necessary and sufficient condition on a Randers space for the existence
of a measure for which Shen's $\bS$-curvature vanishes everywhere.
Moreover, such a measure coincides with the Busemann-Hausdorff measure
up to a constant multiplication.
\end{abstract}

\section{Introduction}\label{sc:intr}

This short article is concerned with a characterization of Randers spaces admitting
measures with vanishing $\bS$-curvature.
A Randers space (due to Randers~\cite{Ra}) is a special kind of Finsler mainfold $(M,F)$
whose Finsler structure $F:TM \lra [0,\infty)$ is written as $F(v)=\alpha(v)+\beta(v)$,
where $\alpha$ is a norm induced from a Riemannian metric on $M$
and $\beta$ is a one-form on $M$.
Randers spaces are important in applications and reasonable for concrete calculations.
See \cite{AIM} and \cite[Chapter~11]{BCS} for more on Randers spaces.

We equip a Finsler manifold $(M,F)$ with an arbitrary smooth measure $m$.
Then the $\bS$-curvature $\bS(v) \in \R$ of $v \in TM$ introduced by Shen
(see \cite[\S 7.3]{Shlec}) measures the difference between $m$ and the volume measure
of the Riemannian structure induced from the tangent vector field of the geodesic $\eta$
with $\dot{\eta}(0)=v$ (see \S \ref{ssc:S} for the precise definition).
The author's recent work \cite{Oint}, \cite{OS} on the weighted Ricci curvature
(in connection with optimal transport theory) shed new light on the importance
of this quantity.

A natural and important question arising from the theory of weighted Ricci curvature
is: when $(M,F)$ admits a measure $m$ with $\bS \equiv 0$?
This is because, if such a measure exists, then we can choose it as a good reference measure.
Our main result provides a complete answer to this question for Randers spaces.

\begin{theorem}\label{th:main}
A Randers space $(M,F)$ admits a measure $m$ with $\bS \equiv 0$
if and only if $\beta$ is a Killing form of constant length.
Moreover, then $m$ coincides with the Busemann-Hausdorff measure
up to a constant multiplication.
\end{theorem}

It has been observed by Shen~\cite[Example~7.3.1]{Shlec} that
a Randers space with the Busemann-Hausdorff measure satisfies $\bS \equiv 0$
if $\beta$ is a Killing form of constant length.
Our theorem asserts that his condition on $\beta$ is also necessary for
the existence of $m$ with $\bS \equiv 0$, and then it immediately follows that
$m$ must be a constant multiplication of the Busemann-Hausdorff measure.

We observe from Theorem~\ref{th:main} that many Finsler manifolds have no
measure with $\bS \equiv 0$.
This means that there is no canonical (reference) measure
(in respect of the weighted Ricci curvature) on such a Finsler manifold.
Therefore, for a general Finsler manifold, it is natural to start with an arbitrary measure,
as was discussed in \cite{Oint} and \cite{OS}.

\section{Preliminaries for Finsler geometry}\label{sc:pre}

We first review the basics of Finsler geometry.
Standard references are \cite{BCS} and \cite{Shlec}.
We will follow the notations in \cite{BCS} with a little change
(e.g., we use $v^i$ instead of $y^i$).

\subsection{Finsler structures}

Let $M$ be a connected $n$-dimensional $C^{\infty}$-manifold with $n \ge 2$,
and $\pi:TM \lra M$ be the natural projection.
Given a local coordinate $(x^i)_{i=1}^n:U \lra \R^n$ on an open set $U \subset M$,
we will always denote by $(x^i;v^i)_{i=1}^n$ the local coordinate on $\pi^{-1}(U)$
given by $v=\sum_i v^i(\del/\del x^i)|_{\pi(v)}$.

A {\it $C^{\infty}$-Finsler structure} is a function $F:TM \lra [0,\infty)$
satisfying the following conditions:
\begin{enumerate}[(I)]
\item $F$ is $C^{\infty}$ on $TM \setminus \{ 0 \}$;
\item $F(cv)=cF(v)$ for all $v \in TM$ and $c \ge 0$;
\item The matrix
\[ g_{ij}(v):=\frac{1}{2} \frac{\del(F^2)}{\del v^i \del v^j}(v) \]
is positive definite for all $v \in TM \setminus \{0\}$.
\end{enumerate}
The positive definite matrix $(g_{ij}(v))$ defines a Riemannian structure $g_v$ of $T_xM$ through
\begin{equation}\label{eq:gv}
g_v\bigg( \sum_i a^i \frac{\del}{\del x^i},\sum_j b^j \frac{\del}{\del x^j} \bigg)
 :=\sum_{i,j} g_{ij}(v) a^i b^j.
\end{equation}
Note that $g_v(v,v)=F(v)^2$.
This inner product $g_v$ is regarded as the best approximation of $F|_{T_xM}$ in the direction $v$.
Indeed, the unit sphere of $g_v$ is tangent to that of $F|_{T_xM}$ at $v/F(v)$
up to the second order.
If $(M,F)$ is Riemannian, then $g_v$ always coincides with the original Riemannian metric.
As usual, $(g^{ij})$ will stand for the inverse matrix of $(g_{ij})$.

Define the {\it Cartan tensor}
\[ A_{ijk}(v):=\frac{F(v)}{2} \frac{\del g_{ij}}{\del v^k}(v) \]
for $v \in TM \setminus \{0\}$, and recall that $A_{ijk} \equiv 0$
if and only if $(M,F)$ is Riemannian.
We also define the {\it formal Christoffel symbol}
\[ \gamma^i{}_{jk}(v) :=\frac{1}{2}\sum_l g^{il}(v) \bigg\{ \frac{\del g_{lj}}{\del x^k}(v)
 +\frac{\del g_{kl}}{\del x^j}(v) -\frac{\del g_{jk}}{\del x^l}(v) \bigg\} \]
for $v \in TM \setminus \{0\}$.
Then the geodesic equation is written as $\ddot{\eta}+G(\dot{\eta})=0$
with the {\it geodesic spray coefficients}
\[ G^i(v):=\sum_{j,k} \gamma^i{}_{jk}(v) v^j v^k \]
for $v \in TM$ ($G^i(0):=0$ by convention).
Using these, we further define the {\it nonlinear connection}
\[ N^i{}_j(v):=\sum_k \bigg\{ \gamma^i{}_{jk}(v)v^k
 -\frac{1}{F(v)} A^i{}_{jk}(v)G^k(v) \bigg\} \]
for $v \in TM$ ($N^i{}_j(0):=0$ by convention),
where $A^i{}_{jk}(v):=\sum_l g^{il}(v)A_{ljk}(v)$.
We can rewrite this as (see \cite[Exercise~2.3.3]{BCS})
\[ N^i{}_j(v)=\frac{1}{2} \frac{\del G^i}{\del v^j}(v). \]

\subsection{$\bS$-curvature and weighted Ricci curvature}\label{ssc:S}

We choose an arbitrary positive $C^{\infty}$-measure $m$ on a Finsler manifold $(M,F)$.
Fix a unit vector $v \in F^{-1}(1)$ and let $\eta:(-\ve,\ve) \lra M$ be the geodesic
with $\dot{\eta}(0)=v$.
Along $\eta$, the tangent vector field $\dot{\eta}$ defines the Riemannian metric
$g_{\dot{\eta}}$ as $(\ref{eq:gv})$.
Denoting the volume form of $g_{\dot{\eta}}$ by $\vol_{\dot{\eta}}$,
we decompose $m$ into $m(dx)=e^{-\Psi(\dot{\eta})} \vol_{\dot{\eta}}(dx)$ along $\eta$.
Then we define the {\it $\bS$-curvature} of $v$ by
\[ \bS(v):=\frac{d(\Psi \circ \dot{\eta})}{dt}(0). \]
We extend this definition to all $w=cv$ with $c \ge 0$ by $\bS(w):=c\bS(v)$.
Clearly $\bS \equiv 0$ holds on Riemannian manifolds with the volume measure.

The {\it weighted Ricci curvature} (introduced in \cite{Oint}) is defined in a similar manner
as follows:
\begin{enumerate}[(i)]
\item $\Ric_n(v):=\Ric(v)+(\Psi \circ \eta)''(0)$ if $\bS(v)=0$,
 $\Ric_n(v):=-\infty$ otherwise;
\item $\Ric_N(v):=\Ric(v)+(\Psi \circ \eta)''(0)-\bS(v)^2/(N-n)$ for $N \in (n,\infty)$;
\item $\Ric_{\infty}(v):=\Ric(v)+(\Psi \circ \eta)''(0)$.
\end{enumerate}
Here $\Ric(v)$ is the usual (unweighted) Ricci curvature of $v$.
The author~\cite{Oint} shows that bounding $\Ric_N$ from below by $K \in \R$
is equivalent to the curvature-dimension condition $\CD(K,N)$,
and then there are rich analytic and geometric applications.
Observe that $\Ric_n \ge K >-\infty$ makes sense only when
the $\bS$-curvature vanishes everywhere.
Therefore the class of such special triples $(M,F,m)$ deserves a particular interest.
We remark that, if two measures $m_1, m_2$ on $(M,F)$ satisfy $\bS \equiv 0$,
then $m_1=c \cdot m_2$ for some positive constant $c$.

We rewrite $\bS(v)$ according to \cite[\S 7.3]{Shlec} for ease of later calculation.
Fix a local coordinate $(x^i)_{i=1}^n$ containing $\eta$ and represent $m$ along $\eta$ as
\[ m(dx)=\sigma(\eta)\, dx^1 dx^2 \cdots dx^n
 =\frac{\sigma(\eta)}{\sqrt{\det(g_{\dot{\eta}})}} \vol_{\dot{\eta}}(dx). \]
Then by definition we have
\[ \bS(v) =\frac{d}{dt}\Big|_{t=0}
 \log\bigg( \frac{\sqrt{\det(g_{\dot{\eta}(t)})}}{\sigma(\eta(t))} \bigg)
 =\frac{1}{2\det(g_v)} \frac{d}{dt}\Big|_{t=0} \big[ \det(g_{\dot{\eta}(t)}) \big]
 -\sum_i \frac{v^i}{\sigma(x)} \frac{\del \sigma}{\del x^i}(x). \]
Since $\eta$ solves the geodesic equation $\ddot{\eta}+G(\dot{\eta})=0$,
the first term is equal to
\begin{align*}
&\frac{1}{2} \sum_{i,j,k} \bigg\{ g^{ij}(v) \frac{\del g_{ij}}{\del x^k}(v) v^k
 +g^{ij}(v) \frac{g_{ij}}{\del v^k}(v) \ddot{\eta}^k(0) \bigg\} \\
&=\sum_{i,k} \bigg\{ \gamma^i{}_{ik}(v) v^k -\frac{1}{F(v)}A^i{}_{ik}(v)G^k(v) \bigg\}
 =\sum_i N^i{}_i(v).
\end{align*}
Thus we obtain
\begin{equation}\label{eq:SN}
\bS(v)=\sum_i \bigg\{ N^i{}_i(v) -\frac{v^i}{\sigma(x)} \frac{\del \sigma}{\del x^i}(x) \bigg\}.
\end{equation}
Note that $\bS(cv)=c\bS(v)$ indeed holds for $c \ge 0$ in this form.

\subsection{Busemann-Hausdorff measure and Berwald spaces}

Different from the Riemannian case, there are several constructive measures
on a Finsler manifold.
Each of them is canonical in some sense, and coincides with
the volume measure on Riemannian manifolds.
Among them, here we treat only the Busemann-Hausdorff measure
which is actually the Hausdorff measure associated with
the suitable distance structure of $(M,F)$.

Using a basis $w_1,w_2,\ldots,w_n \in T_xM$ and its dual basis
$\theta^1,\theta^2,\ldots,\theta^n \in T_x^*M$, the {\it Busemann-Hausdorff measure}
$m_{BH}(dx)=\sigma_{BH}(x)\, \theta^1 \wedge \theta^2 \wedge \cdots \wedge \theta^n$
is defined as
\[ \frac{\omega_n}{\sigma_{BH}(x)} =
 \vol_n \bigg( \bigg\{ (c^i) \in \R^n \,\Big|\, F\bigg( \sum_i c^i w_i \bigg) <1 \bigg\} \bigg), \]
where $\vol_n$ is the Lebesgue measure and $\omega_n$ is
the volume of the unit ball in $\R^n$.

Let $(M,F)$ be a {\it Berwald space} (see \cite[Chapter 10]{BCS} for the precise definition).
Then it is well known that $\bS \equiv 0$ for the Busemann-Hausdorff measure
(see \cite[Proposition~7.3.1]{Shlec}).
In fact, along any geodesic $\eta:[0,l] \lra M$, the parallel transport
$T_{0,t}:T_{\eta(0)}M \lra T_{\eta(t)}M$ with respect to $g_{\dot{\eta}}$
preserves $F$.
Therefore choosing parallel vector fields along $\eta$ as a basis yields that
$\sigma_{BH}$ is constant on $\eta$.
This implies $\bS \equiv 0$.

\section{Proof of Theorem~\ref{th:main}}

Let $(M,F)$ be a Randers space, i.e., $F(v)=\alpha(v)+\beta(v)$ such that
$\alpha$ is a norm induced from a Riemannian metric and that $\beta$ is a one-form.
In a local coordinate $(x^i)_{i=1}^n$, we can write
\[ \alpha(v) =\sqrt{\sum_{i,j} a_{ij}(x) v^i v^j}, \qquad
 \beta(v) =\sum_i b_i(x)v^i \]
for $v \in T_xM$.
The {\it length} of $\beta$ at $x$ is defined by
$\| \beta \|(x):=\sqrt{\sum_{i,j}a^{ij}(x)b_i(x)b_j(x)}$,
which is necessarily less than $1$ in order to guarantee $F>0$ on $TM \setminus \{0\}$.

We denote the Christoffel symbol of $(a_{ij})$ by $\tilde{\gamma}^i{}_{jk}$.
We also define
\[ b^i(x):=\sum_j a^{ij}(x)b_j(x), \qquad b_{i|j}(x):=
 \frac{\del b_i}{\del x^j}(x)-\sum_k b_k(x)\tilde{\gamma}^k{}_{ij}(x). \]
Note that
\begin{equation}\label{eq:beta}
\frac{\del(\| \beta \|^2)}{\del x^i}(x) =2\sum_j b_{j|i}(x) b^j(x).
\end{equation}
We say that $\beta$ is a {\it Killing form} if $b_{i|j}+b_{j|i} \equiv 0$
on $M$.
The geodesic spray coefficients of $F$ are given by (see \cite[(11.3.11)]{BCS})
\begin{align}
G^i(v) &= \sum_{j,k} \gamma^i{}_{jk}(v) v^j v^k \nonumber\\
&= \sum_{j,k} \bigg[ \tilde{\gamma}^i{}_{jk}(x) v^j v^k
 +b_{j|k}(x)\big( a^{ij}(x)v^k -a^{ik}(x)v^j \big) \alpha(v) \nonumber\\
&\quad +b_{j|k}(x) \frac{v^i}{F(v)}
 \big\{ v^j v^k +\big( b^k(x)v^j -b^j(x)v^k \big) \alpha(v) \big\} \bigg] \nonumber\\
&=: \sum_{j,k} \tilde{\gamma}^i{}_{jk}(x) v^j v^k +X^i(v)+Y^i(v). \label{eq:G}
\end{align}

If $\bS \equiv 0$ on $T_xM$, then we deduce from $(\ref{eq:SN})$ that
$\sum_i N^i{}_i(v)$ is linear in $v \in T_xM$.
We shall see that only this infinitesimal constraint is enough to imply
the condition on $\beta$ stated in Theorem \ref{th:main}.
To see this, we calculate $2N^i{}_i=\del G^i/\del v^i$ using $(\ref{eq:G})$.
As the first term $\sum_{j,k}\tilde{\gamma}^i{}_{jk}(x) v^j v^k$ comes from
a Riemannian structure, it suffices to consider only the linearly of
$\sum_i\{ \del X^i/\del v^i(v) +\del Y^i/\del v^i(v) \}$.
For the sake of simplicity, we will omit evaluations at $x$ and $v$
in the following calculations.

We first observe
\[ \sum_i \frac{\del X^i}{\del v^i} =\sum_{i,j} (b_{j|i}-b_{i|j})a^{ij} \alpha
 +\sum_{i,j,k,l}b_{j|k}(a^{ij}v^k -a^{ik}v^j) \frac{a_{il}v^l}{\alpha}
 =0. \]
As Euler's theorem~\cite[Theorem 1.2.1]{BCS} ensures
\[ \sum_i \frac{\del}{\del v^i} \frac{v^i}{F}
 =\frac{1}{F^2} \sum_i \bigg( F- v^i\frac{\del F}{\del v^i} \bigg)=\frac{n-1}{F}, \]
we next obtain
\begin{align*}
\sum_i \frac{\del Y^i}{\del v^i} &=\sum_{i,j} \frac{v^i}{F} \bigg\{
(b_{i|j}+b_{j|i})v^j +(b_{i|j}-b_{j|i})b^j \alpha
 +\sum_{k,l} b_{j|k}(b^k v^j -b^j v^k) \frac{a_{il}v^l}{\alpha} \bigg\} \\
&\qquad +\frac{n-1}{F}\sum_{j,k}b_{j|k}
 \big\{ v^j v^k+(b^k v^j-b^j v^k)\alpha \big\} \\
&=\frac{n+1}{2} \sum_{i,j} (b_{i|j}+b_{j|i}) \frac{v^i v^j}{F}
 +(n+1)\sum_{i,j} (b_{i|j}-b_{j|i})b^j \frac{\alpha v^i}{F}.
\end{align*}
By comparing the evaluations at $v$ and $-v$, the coefficients
$b_{i|j}+b_{j|i}$ in the first term must vanish for all $i,j$,
namely $\beta$ is a Killing form.
For the second term, we find that $(\alpha/F) \sum_j(b_{i|j}-b_{j|i})b^j$
must be constant on each $T_xM$.
If $\alpha/F$ is not constant on some $T_xM$ (i.e., $\| \beta \|(x) \neq 0$), then it holds that
$\sum_j (b_{i|j}-b_{j|i})b^j=0$.
Since $\beta$ is a Killing form, we deduce from $(\ref{eq:beta})$ that
\[ 0=\sum_j(b_{i|j}-b_{j|i})b^j =-2\sum_j b_{j|i}b^j
 =-\frac{\del(\| \beta \|^2)}{\del x^i}. \]
Therefore $\beta$ has a constant length as required,
for $\|\beta\| \neq 0$ is an open condition.
If $\alpha/F$ is constant on some $T_xM$, then the above argument
yields that $\beta \equiv 0$ on $M$.
This completes the proof of the ``only if'' part of Theorem~\ref{th:main}.

For the ``if'' part, it is sufficient to show that the Busemann-Hausdorff measure
satisfies $\bS \equiv 0$, that can be found in \cite[Example~7.3.1]{Shlec}.
We briefly repeat his discussion for completeness.
We first observe from \cite[(2.10)]{Shlec} that
\[ m_{BH}(dx)=\big( 1-\| \beta \|(x)^2 \big)^{(n+1)/2} \sqrt{\det(a_{ij}(x))}\, dx^1 \cdots dx^n
 =: \sigma_{BH}(x)\, dx^1 \cdots dx^n. \]
Since $\beta$ has a constant length, we have
\[ \sum_k \frac{v^k}{\sigma_{BH}(x)} \frac{\del \sigma_{BH}}{\del x^k}(x)
 =\frac{1}{2} \sum_{i,j,k} v^k a^{ij}(x) \frac{\del a_{ij}}{\del x^k}(x)
 =\sum_{i,j} \tilde{\gamma}^i{}_{ij}(x) v^j. \]
Therefore we conclude
\[ \bS(v)=\frac{1}{2} \sum_{i,j,k} \frac{\del}{\del v^i}
 \Big[ \tilde{\gamma}^i{}_{jk}(x) v^j v^k \Big]
 -\sum_k \frac{v^k}{\sigma_{BH}(x)} \frac{\del \sigma_{BH}}{\del x^k}(x) =0. \]
$\qedd$

We finally remark related known results and several consequences of Theorem~\ref{th:main}.

\begin{remark}\label{rm:last}
(a) A Randers space is a Berwald space if and only if $\beta$ is {\it parallel}
in the sense that $b_{i|j} \equiv 0$ for all $i,j$ (see \cite[Theorem 11.5.1]{BCS}).
Thanks to \cite[Example 7.3.2]{Shlec}, we know that a Killing form of constant length
is not necessarily parallel.

(b) In \cite{De}, Deng gives a characterization of vanishing $\bS$-curvature
for homogeneous Randers spaces with the Busemann-Hausdorff measure.

(c) It is easy to construct a Randers space whose $\beta$ does not have a constant length.
Hence many Finsler manifolds do not admit a measure with $\bS \equiv 0$
(in other words, with $\Ric_n \ge K>-\infty$).

(d) Another consequence of Theorem~\ref{th:main} is that only
(constant multiplications of) the Busemann-Hausdorff measures can satisfy
$\bS \equiv 0$ on Randers spaces.
Then a natural question is the following:
\begin{question}
Is there a Finsler manifold $(M,F)$ on which some measure $m$ other than
$($a constant multiplication of$)$ the Busemann-Hausdorff measure satisfies $\bS \equiv 0$?
If yes, what kind of measure is $m$?
\end{question}
If such a measure exists, then it is more natural than the Busemann-Hausdorff
measure in respect of the weighted Ricci curvature.
\end{remark}


{\small
\begin{thebibliography}{AIM}

\bibitem[AIM]{AIM}
P.~L.~Antonelli, R.~S.~Ingarden and M.~Matsumoto,
The theory of sprays and Finsler spaces with applications in physics and biology,
Kluwer Academic Publishers Group, Dordrecht, 1993.

\bibitem[BCS]{BCS}
D.~Bao, S.-S.~Chern and Z.~Shen, An introduction to Riemann-Finsler geometry,
Springer-Verlag, New York, 2000.

\bibitem[De]{De}
S.~Deng, {\it The S-curvature of homogeneous Randers spaces},
Differential Geom.\ Appl.\ {\bf 27} (2009), 75--84.

\bibitem[Oh]{Oint}
S.~Ohta, {\it Finsler interpolation inequalities},
to appear in Calc.\ Var.\ Partial Differential Equations.

\bibitem[OS]{OS}
S.~Ohta and K.-T.~Sturm, {\it Heat flow on Finsler manifolds},
Comm.\ Pure Appl.\ Math.\ {\bf 62} (2009), 1386--1433.

\bibitem[Ra]{Ra}
G.~Randers, {\it On an asymmetrical metric in the fourspace of general relativity},
Phys.\ Rev.\ (2) {\bf 59} (1941), 195--199.

\bibitem[Sh]{Shlec}
Z.~Shen, Lectures on Finsler geometry, World Scientific Publishing Co., Singapore, 2001.

\end{thebibliography}
}
\end{document}